\documentstyle[11pt,epic,eepic]{article} 

\newtheorem{theorem}{Theorem}[section]                                    
\newtheorem{definition}[theorem]{Definition} 
\newtheorem{lemma}[theorem]{Lemma}

\newcommand\bbR{{\bf R}}

\begin{document}     

\title{Stability  of symmetric tops via one variable calculus\\ dg-ga/9608010}

\author{Eugene Lerman 
\thanks{This research was partially supported by
a National Science Foundation Postdoctoral Fellowship.} 
\\ \small Department
of Mathematics, University of Illinois, Urbana, IL 61801 \\ \quad
lerman@math.uiuc.edu
}
 
\maketitle

\centerline{\small AMS classification scheme numbers: 58F}

\begin{abstract}
We study the stability of symmetric trajectories of a particle on the
Lie group $SO(3)$ whose motion is governed by an $SO(3)\times SO(2)$
invariant metric and an $SO(2)\times SO(2)$ invariant potential.  Our
method is to reduce the number of degrees of freedom at {\em singular}
values of the $SO(2)\times SO(2)$ momentum map and study the stability
of the equilibria of the reduced systems as a function of spin.  The
result is an elementary analysis of the fast/slow transition in the
Lagrange and Kirchhoff tops.  More generally, since an $SO(2)\times
SO(2)$ invariant potential on $SO(3)$ can be thought of as ${\bf Z}_2$
invariant function on a circle, we get a condition on the second and
fourth derivatives of the potential at the symmetric points that
guarantees that the corresponding system gains stability as the spin
increases.
\end{abstract}

\section{Introduction}
The motivation for writing this paper is two-fold. The first one is to
see the stable/unstable transition in the Lagrange top without
resorting to the machinery of equivariant singularity theory as in
\cite{G}, \cite{CvdM} let alone the more sophisticated techniques of
\cite{LRSM}.  The idea is to exploit the $SO(2)\times SO(2)$ symmetry
and use reduction instead.  The resulting system can then be analyzed
by looking at the signs of two derivatives of a function of one
variable at one particular point. The second motivation is to explain
why it is that {\sl both} the Lagrange top and the Kirchhoff top [BZ]
undergo the Hamiltonian Hopf bifurcation as their spin varies.

I will show that by carrying out reduction at certain singular values of the 
$SO(2)\times SO(2)$ momentum map one can understand the bifurcation in both
systems (in fact in an infinite family of systems)  as a phenomenon of the 
following sort.  Consider a particle on a line whose behavior is governed
by a 1-parameter family of potentials
$$
	f_\lambda (x) = \frac{1}{6} x^4 - x^2 + \lambda ^2 x^2, 
$$
where $\lambda $ is the parameter.  Then for $|\lambda|<1$ the origin
is a local maximum of the potential, hence an unstable equilibrium of
the system.  Global minima occur at $x= \pm (1-\lambda ^2 )^{1/2}$.
As $|\lambda | $ increases, these stable points move closer to the
origin until they collide at $|\lambda|= 1$; the origin becomes a a
global minimum (hence a stable equilibrium).  In other words, the
reduced system undergoes a generic Hamiltonian ${\bf Z}_2$ invariant
bifurcation --- the figure eight bifurcation (cf \cite{GoSt}).

There is another interpretation of the results of this paper.  Namely, 
the systems that are being studied are three degrees of freedom completely 
integrable systems.  Thus the nondegenerate singularities of the corresponding
Lagrangian foliation are, in Eliasson's classification \cite{E}, of elliptic
or of hyperbolic type.  It turns out that both types of singularities are
present in our systems and that the closures of the sets of elliptic and 
hyperbolic points intersect at some degenerate singular point. This fact
was already observed by Eliasson in the case of the Lagrange top (op.\ cit.).
The point I would like to make is that the Hamiltonian Hopf bifurcation in 
these systems can be described as a transition from an elliptic to a 
hyperbolic singularity through a degenerate singularity.  It would be 
interesting  to see if such transitions are in some sense generic for the 
three degrees of freedom completely integrable systems. In higher dimensions 
many more types of nondegenerate singularities of the Lagrangian
foliation can occur. This suggests that transitions between these singularities
should give rise to interesting bifurcations.

\section{The set-up and the first $SO(2)$ reduction}

The configuration space of the systems we will study is the Lie group $SO(3)$
and the Hamiltonians $H$ are of the form Kinetic Energy + Potential Energy,
where the kinetic term comes from a left $SO(3)$ and right $SO(2)$ invariant
metric, i.e. the moments of inertia satisfy $I_1 = I_2 \not = I_3$, and the 
potential is left-right $SO(2)$ invariant.  Here $SO(2)$ is identified with the
isotropy group of $e_3$, the third vector in the standard basis
$\{e_1, e_2, e_3\}$ of  ${\bf R}^3$.
For example  $W(A) = (Ae_3, e_3)$ is the potential of the Lagrange top and 
$W(A) = (Ae_3, e_1)^2 + (Ae_3, e_2)^2 + c(Ae_3, e_3)^2$, $c$ a positive 
constant, is the potential of the Kirchhoff top.   

Since the left action of the group $SO(2)$ on the cotangent bundle of
$SO(3)$ is free, we may carry out the Marsden-Weinstein-Meyer
reduction.  It amounts to fixing the rate of rotation of the top about
its axis of symmetry.  A routine calculation shows  the reduced spaces,
as manifolds, are the cotangent bundle of the two sphere $S^2$ with
the reduced symplectic structures of the form
$$
	\omega _\lambda = \omega _{T^*S^2} + \lambda \nu,
$$
where $\lambda $ is the value of the $SO(2)$ momentum map, $\omega
_{T^*S^2} $ is the standard symplectic form on the cotangent bundle of
the sphere and $\nu $ is the standard area form on the sphere (pulled
up to the cotangent bundle).  For example, to compute these reduced
spaces one could use Kummer's results [K] or the fact that $SO(3)$ and
the Euclidean group $E(3)$ form a dual pair.  The form $\nu $ is often
referred to as the magnetic form.  This calculation will not be carried
out in the paper.  The reader may simply assume that we are studying a family
of Hamiltonian systems on $(T^*S^2, \omega_\lambda)$,  $\lambda \in \bbR$.

There is one action of $SO(2) $ left on the reduced space.  If we
identify the sphere with the standard round sphere in ${\bf R}^3$ then
the action under the identifications made above is the lift of the
rotation about the $e_3$ axis.  The reduced Hamiltonians are again of
the form Kinetic Energy + Potential Energy, where the kinetic term
comes from a round metric on the sphere and the potential $V$ is an
$SO(2)$ invariant function.  For example the Lagrange top potential is
$V(x) = x_3 |_{S^2}$ and the Kirchhoff top potential is $V(x) =(x_1^2
+ x_2 ^2 + c x_3) |_{S^2} = (1 + (c-1) x_3^2) |_{S^2}$.  Here $(x_1,
x_2, x_3)$ are the standard coordinates on ${\bf R}^3$.

The momentum map $\Phi _\lambda : (T^*S^2, \omega _\lambda ) \to {\bf R}$ is 
the sum of two momentum maps: one coming from the lifted action of $SO(2)$
on $(T^*S^2, \omega _{T^*S^2})$ and the other coming from the action of 
$SO(2)$ on $(S^2, \lambda \nu)$.  We fix the constant in the definition of 
the momentum maps by requiring the map to be $0$ at the North pole.

We think of the reduced systems $(T^*S^2, \omega _\lambda , h)$ as two
degrees of freedom systems depending on one real parameter $\lambda$,
the spin.  The action of $SO(2)$ fixes the North and South poles.
Consequently these points are critical for our Hamiltonians $h$.  We
would like to understand how the behavior of the systems near the
symmetric points changes as we vary the parameter $\lambda$.  We will
consider the North pole.  The arguments regarding the South pole are
exactly the same.  In fact the Hamiltonian for the Kirchhoff top is
invariant under the ${\bf Z}_2$ symmetry interchanging the poles.

To understand the behavior of the system near the North pole we reduce with 
respect to the remaining symmetry at $\Phi _\lambda=0$, the value of the 
normalized momentum map at the North pole.  Since the North pole is a fixed 
point, the standard Marsden-Weinstein-Meyer procedure does not apply.  
We carry out the reduction in the next section.

\section{Reduction at a critical value of the momentum map}

Let us briefly summarize the procedure that we are going to follow. A
full account is given in \cite{SL}.  Our point of view is that given a
Hamiltonian action of a compact Lie group $K$ on a symplectic manifold
$(M, \omega )$ with momentum map $F:M \to {\bf k}^*$, the reduced
space at a coadjoint orbit $O_\alpha = K\cdot \alpha \subset {\bf
k}^*$ is the Hausdorff topological space
$$
		M_\alpha := F^{-1} (O_\alpha )/K.
$$
Smooth $K$ invariant functions on $M$ descend to continuous functions
on $M_\alpha$; we denote these functions by $C^\infty (M_\alpha )$.
Algebraically
$$
	C^\infty (M_\alpha ) = C^\infty (M)^K/I_\alpha,
$$
where $ C^\infty (M)^K$ denotes the algebra of smooth $K$ invariant
functions and $I_\alpha$ is the ideal in $C^\infty (M)^K$ consisting
of functions that vanish on $F^{-1} (O_\alpha )$.  It is not hard to
see that $I_\alpha $ is a {\sl Poisson } ideal \cite{ACG}.  Indeed,
the Hamiltonian flows of $K$ invariant functions preserve the set
$F^{-1} (O_\alpha )$, so for any $f\in C^\infty (M)^K$ and any $h\in
I_\alpha $ we have $\{f,h\} = 0$.  Therefore $C^\infty (M_\alpha )$ is
a Poisson algebra.  The Poisson bracket allows us to associate flows
to functions.  It also makes it possible for us to say when two
reduced spaces are isomorphic.
\begin{definition}
Two reduced spaces $X$ and $Y$ are {\em isomorphic} if there exists a 
homeomorphism $\phi:X \to Y$ such that the induced map
 $\phi^* : C^\infty (Y) \to C^\infty (X)$ is a Poisson isomorphism.
\end{definition}

The following example will be very important for us.  In fact, thanks
to the equivariant Darboux theorem, this is precisely the situation
that we are interested in. Consider the standard lifted action of
$SO(2)$ on the cotangent bundle $T^* {\bf R}^2$.  Let $(x,y)$ be the
coordinates on ${\bf R}^2$ and $(x,y, p_x, p_y)$ the canonical
coordinates on $T^* {\bf R}^2$.  The standard symplectic form $\omega
$ is given in these coordinates by $\omega = dx\wedge dp_x + dy \wedge
dp_y$ and the momentum map $F$ by
$$
		F(x,y,p_x, p_y) = xp_y - y p_x .
$$
The map $\psi: T^* {\bf R} \to T^*\bbR ^2$, $(u,p_u) \mapsto (u,0,
p_u, 0)$ is a symplectic embedding whose image lies entirely in the
zero level set $F^{-1}(0)$ of $F$.  Moreover, all the orbits of
$SO(2)$ in $F^{-1}(0)$ (except $\{0\}$) intersect the image in exactly
two points: $(u,0, p_u,0)$ and $(-u,0, -p_u,0)$.  Consequently $\psi$
descends to a homeomorphism $\phi : T^* \bbR/{\bf Z}_2 \to
F^{-1}(0)/SO(2) \equiv (T^* \bbR^2)_0$ where ${\bf Z}_2$ acts on $T^*
\bbR$ by the map $(u,p_u) \mapsto (-u,- p_u)$.

Note that since any $SO(2)$ invariant function on $T^* \bbR ^2$ pulls
back via $\psi $ to a ${\bf Z}_2$ invariant function on $T^* \bbR$,
$\phi$ pulls back smooth functions on the reduced space $(T^*
\bbR^2)_0$ to smooth functions on the orbifold $T^* \bbR/{\bf Z}_2$.
Since $\psi $ is a symplectic embedding the pull-back $\phi^*$ is
Poisson.  It is easy to see that $\phi ^*$ is injective: $\phi$ is a
homeomorphism.  To prove surjectivity of the pull-back we use a
theorem of G.\ Schwarz on invariant smooth functions [Sch].  The ${\bf
Z}_2$ invariant polynomials on $T^* \bbR $ are generated by three
polynomials: $u^2$, $p_u ^2$ and $up_u$.  Schwarz's theorem in this
case says that any smooth ${\bf Z}_2 $ invariant function is a
composition of a smooth function on $\bbR^3$ with the map $T^* \bbR
\to \bbR^3$, $(u, p_u) \mapsto (u^2, p_u ^2, up_u)$, i.e., if $f\in
C^\infty (T^* \bbR)^{{\scriptsize \bf Z}_2}\equiv C^\infty (T^*
\bbR/{\bf Z}_2)$ then there is $\bar{f} \in C^\infty (\bbR ^3)$ with
$$
		f(u, p_u) = \bar{f} (u^2, p_u ^2, up_u).
$$ 
Now $u^2 = \psi ^* (x^2 + y^2)$, $p_u ^2 = \psi ^* (p_x ^2 + p_ y^2)$
and $up_u = \psi ^* (xp_x + y p_y)$.  Therefore any function $f\in
C^\infty (T^* \bbR/{\bf Z}_2)$ is in the image of $\phi ^*$.\\ {\sc
Conclusion: } the reduced space $(T^* \bbR^2)_0$ is isomorphic to the
orbifold $ T^* \bbR/{\bf Z}_2$ (the cotangent bundle of a line reduced
by a finite group ${\bf Z}_2$).

\subsection{Stability}

Marsden and Weinstein \cite{MW} proved that stable equilibria of reduced 
systems correspond to stable relative equilibria of the original system.  More 
precisely they considered a Lie group $G$ acting on a symplectic manifold 
$(M, \omega )$ in a  Hamiltonian 
fashion with corresponding momentum map $\Phi : M \to {\bf g}^*$ and an 
invariant Hamiltonian $h\in C^\infty (M)^G$.  For a regular value $\alpha $ of 
the momentum map let $h_\alpha $ denote the corresponding reduced Hamiltonian 
on the reduced space $M_\alpha := \Phi ^{-1}(G\cdot\alpha )/ G$.  
Marsden and Weinstein showed that if $x\in M_\alpha $ is a critical point of 
$h_\alpha$ and  the Hessian of $h_\alpha $
at $x$ is  definite then the trajectory of the original Hamiltonian $h$ through
a point $m\in \Phi^{-1}(\alpha )$ that corresponds to  $x\in M_\alpha $ is 
stable relative to the group action. In fact the result holds under a 
slightly weaker hypothesis, namely that the action of the group $G$ is free at 
$m$, i.e., the reduced space $M_\alpha$ is smooth near $x$.  
Thus $\alpha $ being a regular value is not really needed.
 
If a point $m$ is symmetric, i.e., has a nontrivial isotropy, then it is 
automatically a critical point of any invariant Hamiltonian. 
Since it is also a critical point of the momentum map Marsden and Weinstein's
result does not apply.  However, for the action of $SO(2)$ on the cotangent 
bundle of $\bbR^2$ that we are considering, the correspondence still
holds. 

\begin{lemma}
Suppose a function $h\in C^\infty(T^* \bbR^2)$ is $SO(2)$ invariant.
Then if the Hessian of $h_0 := h|_{\{(u,0,p_u, 0)\in T^* {\scriptsize
\bf R}^2\}}$ is definite at the origin, the origin $0=(0,0,0,0)\in T^*
\bbR^2$ is $SO(2)$ stable, i.e., given any $SO(2)$ invariant
neighborhood $U$ of $0$ in the cotangent bundle $T^* \bbR^2$ there is
an $SO(2)$ invariant neighborhood $V$ of $0$ such that for all $x\in
V$ the trajectory of the Hamiltonian flow of $h$ through $x$ stays in
$U$.
\end{lemma}

\noindent {\sc Proof. }  It is no loss of generality to assume that
the Hessian ${\cal H}(h_0)$ of $h_0$ at $0$ is of the form
$$
		{\cal H}(h_0) = au^2 + b p_u ^2; \quad a,b>0.
$$
Indeed, any quadratic form on $\bbR^2$ can be diagonalized by an
orthogonal change of coordinates, 
$\left(\begin{array}{l}
 u \\
p_u
\end{array} \right) 
\mapsto 
\left( \begin{array}{rr}
c &s\\
-s&c 
\end{array} \right) 
\left(\begin{array}{l}
u\\
p_u
\end{array} \right)$, where $c^2 + s^2 = 1$.  Moreover, this change of
coordinates is induced by an $SO(2)$ equivariant  change of
coordinates on $T^* \bbR^2$; it is given by the matrix
$$
\left( \begin{array}{rrrr}
c &0&s&0\\
0&c&0&s\\
-s&0&c &0\\
0&-s&0& c
\end{array} \right) .
$$
Since the Hessian of $h$ at $0$ is $SO(2)$ invariant and since 
${\cal H} (h) |_{\{(u,0,p_u, 0)\in  T^* \bbR^2\}}={\cal H} (h_0)$, the 
Hessian of $h$ has to be of the form
$$
	{\cal H} (h) = a (x^2 + y^2) + b (p_x ^2 + p_y ^2) + c F,
$$
where $c$ is a constant and $F = xp_y - yp_x$ is the momentum map.  It
cannot have any more terms for then the Hessian of the reduced
Hamiltonian ${\cal H}(h_0)$ would have more terms.  

Let $\rho$ be a a compactly supported smooth function on the cotangent
bundle of $\bbR^2$ which is identically 1 near $0$.  Then near the
origin the Poisson bracket $\{h, \rho cF\}$ is zero.  The Hessian of
$\bar{h} := h - \rho cF$ is positive definite at $0$, hence the origin
is a stable equilibrium of $\bar{h}$.  On the other hand, since the
flows of $\bar{h}$ and of $\rho cF$ commute near the origin, the flow
of the original Hamiltonian is the composition of the flows of
$\bar{h}$ and of $cF$.  Therefore $0 \in T^* \bbR^2$ is $SO(2)$ stable
for the flow of $h$. \hfill $\Box$

\subsection{Reduction of the cotangent bundle of the 2-sphere}

Consider  coordinates on the upper hemisphere of 
$S^2$ induced by the projection $(x,y,z)\mapsto (x,y)$. The inverse
map $\phi $ is
given by  $\phi(x,y) =
\mapsto (x,y,\sqrt{1-x^2-y^2})$.  These coordinates allow us to 
compute the most interesting part of the reduced system.  We don't need the 
coordinates to compute the reduced space per se, rather the coordinates
are convenient for computing the reduced Hamiltonians. Since
$$
\phi ^*(xdy\wedge dz + ydz\wedge dx + z dx \wedge dy) 
		= \frac{1}{z} dx \wedge dy,
$$
where $z= \sqrt{1-x^2-y^2}$, the symplectic form $\omega _\lambda$ in the 
canonical coordinates corresponding to $(x,y)$ is 
$$
\omega _\lambda =dx\wedge dp_x + dy\wedge dp_y +\frac{\lambda }{z} dx\wedge dy.
$$
The momentum map $\Phi _\lambda $ in these coordinates is then 
$$
\Phi _\lambda (x,y,p_x, p_y) = yp_x - xp_y - \lambda z + \lambda.
$$
Note that the North pole in these coordinates is the origin and that 
$\Phi (0) = 0$.

The metric $g$ by assumption is the one that comes from the embedding of 
$S^2$ as a round sphere in $\bbR ^3$.  We now normalize the metric by setting
the radius of the sphere to 1. Then the sphere is cut out by the
equation 
\begin{equation} \label{eq1}
x^2 + y^2+ z^2 = 1,
\end{equation}
 In coordinates the metric $g$ is given by 
$$
g=\left( \begin{array}{lr}
			1+\frac{x^2}{z^2} & \frac{xy}{z^2}\\
			\frac{xy}{z^2}   & 1+\frac{y^2}{z^2}
							\end{array}\right).
$$
Consequently
$$
g^{-1}= z^2
			\left( \begin{array}{lr}
			1+\frac{y^2}{z^2} & -\frac{xy}{z^2}\\
			-\frac{xy}{z^2}   & 1+\frac{x^2}{z^2}
							\end{array}\right).
$$

We are now in position to compute the reduced Hamiltonian.  From the
discussion in the beginning of the section we know that locally the
reduced space is the cotangent bundle of the interval $(-1, 1)$
divided by the involution $(u, p_u) \mapsto (-u, -p_u)$.  Since our
coordinates are not Darboux we proceed as follows.  Set $y=0$.  Then
$\Phi_\lambda (x,0,p_x,p_y) =0$ if and only if
$$
	p_y =\lambda (1-z)/x .
$$
Note that for $y=0$, 
$$
 1-z= 1-\sqrt{1-x^2} = 1-(1-\frac{1}{2}x^2 -\frac{1}{8}x^4 +{\cal O}(x^6))
         = \frac{1}{2}x^2 +\frac{1}{8}x^4 +{\cal O}(x^6) ,
$$
so 
$$
	p_y (x) = \frac{\lambda}{2}(x +\frac{1}{8}x^3 +{\cal O}(x^5))
$$
is smooth at $x=0$.  Consider now a map $\psi:T^* (-1,1) \to 
\Phi_\lambda ^{-1} (0)$, $(u,p_u) \mapsto (u,0,p_u, \lambda (1-z)/u)$ where now
$z=\sqrt{1-u^2}$.  As in the example considered above, the composition
$$
  T^* (-1,1) \to \Phi_\lambda ^{-1} (0) \to \Phi_\lambda ^{-1} (0)/SO(2)
$$
is a branched double cover.  Therefore, we will simply work on $T^* (-1,1)$
keeping the involution in mind.  Now the pull-back by  $\psi$ of the kinetic
energy part of the Hamiltonian is 
$$
\frac{	z^2}{2} (p_u ^2 + \left(1+\frac{x^2}{z^2}\right)
	\left(\lambda \frac{1-z}{u}\right)^2) =
\frac{1}{2}	G(u)p_u^2 + \frac{\lambda^2}{2}\left( \frac{1-z}{u}\right)^2,
$$
where $G(u) = (1-u^2)$.  The potential energy part pulls back to
an even function of $u$, hence to a function of $u^2$.  For example, for the 
Lagrange top the corresponding function is 
$$
	V(u^2) = \sqrt{1-u^2}
$$
and for the Kirchhoff top 
$$
V(u^2) = 1+(c-1)(1-u^2).
$$
The critical points of the reduced Hamiltonian 
$$
	H_\lambda (u,p_u) = \frac{1}{2}G(u)p_u^2 + \frac{1}{2}\lambda ^2
\left(\frac{1-\sqrt{1-u^2}}{u} \right)^2 + V(u^2)
$$
are in one-to-one correspondence with the critical points of the effective 
potential
$$
	U_\lambda (u) 
		= \frac{1}{2}\lambda ^2\left(\frac{1-\sqrt{1-u^2}}{u}\right)^2 + V(u^2).
$$
Moreover, relative maxima of $U_\lambda$ correspond to unstable equilibria
of $H_\lambda$ and relative minima to stable equilibria.\\[3pt]

\noindent
{\sc Remark  } Since we are working on the branched double cover of the 
reduced system, any pair of critical points of $U_\lambda$ of the form
$\pm u$, $u\neq 0$ correspond to the same critical point of the reduced system.
Therefore only the non-negative critical points of the potential need be 
considered.

\subsection{Critical values of the effective potential $U_\lambda (u)$, 
$u\geq 0$}

Since
$$
\left(\frac{1-\sqrt{1-u^2}}{u}\right)^2 = \frac{u^2}{4} (1+ \frac{1}{2} u^2
		+{\cal O}(u^4))
$$
there exists a change of variables $u= \tau (v)$ so that
$$
	v^2 =  \left(\frac{1-\sqrt{1-\tau(v)^2}}{\tau(v)}\right)^2.
$$
Note that $\tau (v)$ is an odd function of the form
$$
	\tau (v)= 2v - 2v^3 + {\cal O}(v^5).
$$
Then the transformed reduced potential $V(\tau(v)^2)$ is even in $v$ and so 
has to be of the form $V(\tau(v)^2) = f(v^2)$ for some smooth function $f$.
Therefore the study of the bifurcation of the critical points of the potential 
$U_\lambda (u)$ for $u$ small and non-negative is reduced to studying the bifurcation of the
critical points of $F_\lambda (v) = \frac{\lambda^2}{2}v^2 + f(v^2)$ for $v$ 
small and nonnegative.

It will be useful to express the derivatives $f'(0)$ and $f''(0)$ in terms
of the derivatives $V'(0)$ and $V''(0)$.  Since $\tau (v)^2 = 4v^2 - 8v^4 +
{\cal O}(v^6)$ and $\tau (v)^4 = 16v^4 + {\cal O}(v^6)$, we have
\begin{eqnarray*}
V(\tau(v)^2) &=& V(0) + V'(0)\tau (v)^2 + V''(0)\tau (v) ^4 + {\cal O}(v^6) \\
	     &=& V(0) + 4V'(0)v^2 +8(V''(0) - V'(0))v^4 +{\cal O}(v^6).
\end{eqnarray*}
Therefore
\begin{equation}\label{eq2}
f'(0) = 4V'(0) \quad \mbox{ and } \quad f'' (0) = 8(V''(0) - V'(0)).
\end{equation}

\subsection{Bifurcation of critical points of 
	$F_\lambda(v) = \frac{\lambda^2}{2} v^2 + f(v^2)$, $v\geq 0$}

It is no loss of generality to assume that $\lambda >0$.  We also make
the genricity assumptions: $f'(0)\not = 0$ and $f''(0)\not = 0$.
$$
	\frac{d}{dv}F_\lambda (v) = v(\lambda ^2 + 2f'(v^2))
$$
Thus if $f'(0)>0$ then $v=0$ is the only critical point for all values of
$\lambda$, i.e., no bifurcations occur. Physically this mean that if a top
in a straight up position with zero spin is stable, then it is stable in this
position for all values of spin.  

Assume now $f'(0)< 0$.  Then 
$f'(0) =- a^2/2$
for some $a>0$.  Consider $\phi (t, \lambda) = \lambda ^2 + 2f'(t)$.
Since $(f')'(0) \not = 0$, by the
implicit function theorem there exists a smooth function $g(\lambda)$ defined
near $\lambda = a$ such that $g(a) = 0 $ and 
$\phi(g(\lambda), \lambda) = 0$, i.e.,
\begin{equation}\label{eq3}
	\frac{1}{2} \lambda ^2 + f'(g(\lambda )) =0.
\end{equation}
Note that since
$$
0 = \left.\frac{d}{d\lambda}\right|_{\lambda =a} \phi (g(\lambda ), \lambda)=
	\left.[\lambda + f''(g(\lambda))g'(\lambda)]\right|_{\lambda =a}
$$
we have 
$$
 f''(0)g'(a) = -a <0 .
$$
Thus if $f''(0)>0$ then $g(\lambda )$ is strictly decreasing near
$\lambda =a$.  Consequently 
there are $\epsilon, \delta >0$ so that for $0\leq v<\delta$  the equation
 $\lambda ^2 + f'(v^2)=0$
has no solution whenever $a<\lambda < a+ \epsilon$ and for 
$a-\epsilon <\lambda < a$ the equation has exactly one solution, namely
$(v,\lambda )= (g(\lambda )^{1/2}, \lambda)$; see figure~1.

\begin{figure}[h]\label{fig}
\setlength{\unitlength}{0.009125in}
\begin{picture}(560,575)(0,-10)
\thicklines
\drawline(393,125)(393,125)
\drawline(392,117)(392,117)
\drawline(389,110)(389,110)
\drawline(387,103)(387,103)
\drawline(382,96)(382,96)
\drawline(372,86)(372,86)
\drawline(364,83)(364,83)
\drawline(356,81)(356,81)
\drawline(394,133)(394,133)
\drawline(395,140)(395,140)
\thinlines
\path(60,460)(60,400)
\path(0,420)(270,420)
\path(262.000,418.000)(270.000,420.000)(262.000,422.000)
\put(350.833,319.167){\arc{78.352}{0.5112}{1.5921}}
\put(227.500,377.500){\arc{351.070}{0.1862}{0.4573}}
\put(86.591,673.864){\arc{654.118}{0.4217}{0.8188}}
\put(-197.500,272.500){\arc{674.129}{5.4138}{5.7631}}
\path(95.000,440.000)(92.422,444.423)(89.778,448.807)
	(87.067,453.150)(84.290,457.451)(81.448,461.710)
	(78.542,465.925)(75.573,470.095)(72.540,474.220)
	(69.445,478.298)(66.289,482.328)(63.071,486.311)
	(59.794,490.244)(56.457,494.126)(53.062,497.958)
	(49.608,501.738)(46.098,505.464)(42.532,509.137)
	(38.910,512.756)(35.234,516.319)(31.504,519.825)
	(27.721,523.275)(23.886,526.667)(20.000,530.000)
\path(0,220)(270,220)
\path(262.000,218.000)(270.000,220.000)(262.000,222.000)
\thicklines
\put(364,83){\blacken\ellipse{2}{2}}
\put(364,83){\ellipse{2}{2}}
\put(372,86){\blacken\ellipse{2}{2}}
\put(372,86){\ellipse{2}{2}}
\put(377,90){\blacken\ellipse{2}{2}}
\put(377,90){\ellipse{2}{2}}
\put(383,96){\blacken\ellipse{2}{2}}
\put(383,96){\ellipse{2}{2}}
\put(387,103){\blacken\ellipse{2}{2}}
\put(387,103){\ellipse{2}{2}}
\put(389,110){\blacken\ellipse{2}{2}}
\put(389,110){\ellipse{2}{2}}
\thinlines
\path(60,460)(60,460)(60,560)
\path(62.000,552.000)(60.000,560.000)(58.000,552.000)
\thicklines
\put(356,81){\blacken\ellipse{2}{2}}
\put(356,81){\ellipse{2}{2}}
\put(394,133){\blacken\ellipse{2}{2}}
\put(394,133){\ellipse{2}{2}}
\put(392,117){\blacken\ellipse{2}{2}}
\put(392,117){\ellipse{2}{2}}
\dashline{6.000}(350,80)(350,20)
\thinlines
\path(290,20)(560,20)
\path(552.000,18.000)(560.000,20.000)(552.000,22.000)
\thicklines
\path(60,90)(60,150)
\dashline{6.000}(60,90)(60,25)
\thinlines
\path(0,20)(270,20)
\path(262.000,18.000)(270.000,20.000)(262.000,22.000)
\path(60,150)(60,165)
\path(62.000,157.000)(60.000,165.000)(58.000,157.000)
\path(60,20)(60,0)
\path(350,20)(350,0)
\path(350,145)(350,165)
\path(352.000,157.000)(350.000,165.000)(348.000,157.000)
\thicklines
\path(350,80)(350,145)
\drawline(375,90)(375,90)
\drawline(375,90)(375,90)
\thinlines
\path(60,260)(60,200)
\path(60,260)(60,260)(60,360)
\path(62.000,352.000)(60.000,360.000)(58.000,352.000)
\path(290,220)(560,220)
\path(552.000,218.000)(560.000,220.000)(552.000,222.000)
\path(350,260)(350,200)
\path(350,260)(350,260)(350,360)
\path(352.000,352.000)(350.000,360.000)(348.000,352.000)
\path(290,420)(560,420)
\path(552.000,418.000)(560.000,420.000)(552.000,422.000)
\path(350,460)(350,400)
\thicklines
\put(393,125){\blacken\ellipse{2}{2}}
\put(393,125){\ellipse{2}{2}}
\thinlines
\path(350,460)(350,460)(350,560)
\path(352.000,552.000)(350.000,560.000)(348.000,552.000)
\thicklines
\put(395,143){\blacken\ellipse{2}{2}}
\put(395,143){\ellipse{2}{2}}
\thinlines
\path(60,290)	(63.243,290.157)
	(65.000,290.000)

\path(65,290)	(67.391,289.241)
	(70.231,288.025)
	(75.000,285.000)

\path(75,285)	(77.923,280.150)
	(79.125,277.300)
	(80.000,275.000)

\path(80,275)	(81.318,270.631)
	(82.008,267.910)
	(82.693,265.033)
	(83.353,262.152)
	(83.971,259.419)
	(85.000,255.000)

\path(85,255)	(85.734,252.062)
	(86.734,248.062)
	(87.370,245.517)
	(88.117,242.531)
	(88.989,239.045)
	(90.000,235.000)

\thicklines
\path(60,90)	(63.243,90.157)
	(65.000,90.000)

\path(65,90)	(67.391,89.241)
	(70.231,88.025)
	(75.000,85.000)

\path(75,85)	(77.923,80.150)
	(79.125,77.300)
	(80.000,75.000)

\path(80,75)	(81.318,70.631)
	(82.008,67.910)
	(82.693,65.033)
	(83.353,62.152)
	(83.971,59.419)
	(85.000,55.000)

\path(85,55)	(85.734,52.062)
	(86.734,48.062)
	(87.370,45.517)
	(88.117,42.531)
	(88.989,39.045)
	(90.000,35.000)

\put(395,590){\makebox(0,0)[lb]{$f''(0)<0$}}
\put(335,480){\makebox(0,0)[lb]{$a$}}
\put(330,280){\makebox(0,0)[lb]{$a$}}
\put(330,85){\makebox(0,0)[lb]{$a$}}
\put(380,485){\makebox(0,0)[lb]{$g(\lambda ) = v^2$}}
\put(390,285){\makebox(0,0)[lb]{$g(\lambda )^{1/2} = v$}}%
\put(360,550){\makebox(0,0)[lb]{$\lambda $}}
\put(360,155){\makebox(0,0)[lb]{$\lambda $}}
\put(350,350){\makebox(0,0)[lb]{$\lambda $}}
\put(550,205){\makebox(0,0)[lb]{$v$}}
\put(550,5){\makebox(0,0)[lb]{$v$}}
\put(550,400){\makebox(0,0)[lb]{$v^2$}}

\put(70,550){\makebox(0,0)[lb]{$\lambda$}}
\put(40,290){\makebox(0,0)[lb]{$a$}}
\put(40,85){\makebox(0,0)[lb]{$a$}}
\put(40,480){\makebox(0,0)[lb]{$a$}}
\put(90,460){\makebox(0,0)[lb]{$g(\lambda) = v^2$}}
\put(95,275){\makebox(0,0)[lb]{$g(\lambda)^{1/2} = v$}}

\put(70,155){\makebox(0,0)[lb]{$\lambda$}}
\put(65,350){\makebox(0,0)[lb]{$\lambda$}}
\put(260,400){\makebox(0,0)[lb]{$v^2$}}
\put(255,205){\makebox(0,0)[lb]{$v$}}
\put(260,5){\makebox(0,0)[lb]{$v$}}
\put(105,590){\makebox(0,0)[lb]{$f''(0)>0$}}%
\end{picture}
\begin{center}
Figure~1: bifurcation of critical points of $F_\lambda (v)$
\end{center}
\end{figure}

Similarly if  $f''(0)<0$ then $g(\lambda )$ is strictly increasing near
$\lambda =a$.  Consequently there are $\epsilon, \delta >0$ so that for 
$0\leq v<\delta$  the equation  $\lambda ^2 + f'(v^2)=0$ has no solution 
whenever $a -\epsilon <\lambda < a$ and for 
$a <\lambda < a+\epsilon $ the equation has exactly one solution, namely
$(v,\lambda )= (g(\lambda )^{1/2}, \lambda).$

Finally we determine the stability of these critical points.  Since
$$
	\frac{d^2}{dv^2} F_\lambda (v) = \lambda ^2 + 2f'(v^2) 
		+ 4v^2 f''(v^2)
$$
we have, using (3),
$$
	\left.\frac{d^2}{dv^2} \right|_{v= g(\lambda )^{1/2}} F_\lambda (v) =
		4|g(\lambda )| f''(g(\lambda ))\quad \mbox{which is}\quad
\left\{\begin{array}{c}
	>0\, \mbox{ if } f''(0)>0\\
	<0 \, \mbox{ if } f''(0)<0\\
	\end{array}\right.
$$
and
$$
	\left.\frac{d^2}{dv^2} \right|_{v= 0} F_\lambda (v) =
		2(\lambda ^2 + f'(0))
\quad \mbox{which is}\quad
\left\{\begin{array}{c}
	>0 \,\mbox{ if } \lambda > \sqrt{-2f'(0)}\\
	<0 \,\mbox{ if } \lambda <\sqrt{-2f'(0)} \\
	\end{array}\right.
$$
\\[2pt]
\noindent
{\sc Conclusions } Assume that $f'(0)<0$.\\
If $f''(0)> 0$ then as $\lambda $ decreases below $\lambda_0 = \sqrt{-2f'(0)}$,
$v=0$ changes from a local minimum to a local 
maximum and a one parameter family of local minima 
$v = g(\lambda )^{1/2}$ bifurcate from $v=0$.

If $f''(0)<0$ then the situation is reversed.  For
 $\lambda < \lambda_0$ and $v$ small,
the origin is a unique critical point, and it is a local maximum.  As $\lambda$
increases past the critical value, $v=0$ becomes a local minimum and a one
parameter family $v= g(\lambda )^{1/2}$ of local maxima bifurcate from it.\\

These results easily translate into a description of the bifurcation of 
nonnegative critical points of the potential
$$
U_\lambda (u) = 
  \frac{\lambda ^2}{2}\left( \frac{1 -\sqrt{1-u^2}}{u} \right)^2 + V(u^2),
$$
where without loss of generality we shall assume that $\lambda \geq 0$.
Recall that
$$
f'(0) = 4V'(0) \quad \mbox{ and } \quad f'' (0) = 8(V''(0) - V'(0)).
$$
It follows that if $V'(0)>0$ then $u=0$ is a local minimum for all values of 
$\lambda$ and no bifurcations occur. We now assume that $V'(0)<0$. Then there
are two alternatives.\\[4pt]
1. If $V''(0)> V'(0)$ then as $\lambda $ decreases below 
$\lambda _0 = \sqrt{-V'(0)} $, $u=0$ changes from a local minimum to a local
maximum of the potential $U_\lambda $ and a 1-parameter family 
$u = \rho (\lambda )$ of local minima bifurcate off from $u=0$.  Here 
$\rho (\lambda)$ is smooth for $\lambda <\lambda _0$ and H\"{o}lder $1/2$
at $\lambda_0$.  Consequently the corresponding system on the cotangent bundle
of the two-sphere undergoes a Hamiltonian Hopf bifurcation.\\[4pt]
2. If $V''(0)< V'(0)$ then for $0\leq \lambda < \lambda _0$ and $u$ small, the 
point $u=0$ is a unique local maximum.  As $\lambda $ increases past 
$\lambda _0$, $u=0$ becomes a local minimum and a 1-parameter family 
$u=\rho (\lambda )$ of local maxima bifurcate from it.\\[4pt]
For the Lagrange top $V(u^2) = \sqrt{1-u^2}$, so the first alternative
holds with $\lambda _0 = 2$ (this agrees with \cite{G}).  For the Kirchhoff
top $V(u^2)= 1+ (c-1)(1-u^2)$.  So if $c<1$ there is no bifurcation and if
$c>1$ the first alternative holds.

\subsection*{Acknowledgments}

I would like to thank Victor Guillemin for getting me interested in tops and
for making available to me his unpublished notes on the Lagrange top.  I 
would also like to thank Iliya Zakharevich for his helpful comments and 
Richard Cushman for the encouragement.

\end{document}